\documentclass[12pt]{article}
\usepackage{amssymb,amsfonts,amsmath,amsthm,cite}
\usepackage{epsfig}
\newtheorem{thm}{Theorem}[section]
\newtheorem{cor}{Corollary}[section]
\newtheorem{lem}{Lemma}[section]

\makeatletter \@addtoreset{equation}{section}

\def\qed{\hfill \rule{4pt}{7pt}}

\begin{document}
{\Large
\begin{center}
 The Ratio  Monotonicity  of the  $q$-Derangement Numbers
\end{center}
}

\begin{center}
William Y. C. Chen$^1$ and Ernest X. W. Xia$^2$ \\
Center for Combinatorics, LPMC-TJKLC\\
Nankai University\\
 Tianjin 300071, P. R. China\\
Email: $^1$chen@nankai.edu.cn,  $^2$xia@cfc.nankai.edu.cn

\end{center}


\noindent {\bf Abstract.} We show that the  $q$-derangement number
$d_n(q)$ $(n \geq 6)$  satisfies a ratio monotone  property,  which
implies the spiral property and the log-concavity except for the
last term when $n$ is
 even. More specifically, we discover two increasing ratio sequences
motivated by the spiral property.



\section{Introduction}

Let $D_n$ be the set of derangements on $\{1,2,\ldots,n\}$, and let
$\rm{maj}(\pi)$ denote the major index of a permutation $\pi$. The
$q$-derangement number $d_n(q)$ is defined as
\[ \sum\limits_{\pi
\in D_n}q^{\rm{maj}(\pi)}.
\]
 The following formula is due to Wachs \cite{Wachs} (see also Gessel and
 Reutenauer \cite{GR}):
\[
d_n(q)=[n]!\sum_{k=0}^n(-1)^kq^{k\choose 2}\frac{1}{[k]!},
\]
where $[n]=1+q+q^2+\cdots +q^{n-1}$ and $[n]!=[1] [2] \cdots [n]$.
Note that $d_n(q)$ satisfies the following recursion:
\begin{align}\label{recursion}
d_n(q)=(1+q+q^2+\cdots +q^{n-1})d_{n-1}(q)+(-1)^nq^{n\choose 2},
 \ \ \ \ \ \ \ \ n\geq 2,
 \end{align}
with the initial value $d_1(q)=0$.

Chen and Rota \cite{Chen} showed that the $q$-derangement numbers
are unimodal and conjectured the maximum coefficient appears in the
middle. Zhang \cite{Zhang} confirmed this conjecture by showing that
the $q$-derangement numbers satisfy the spiral property. For
example, we have
 \begin{align*}
d_8(q)=q&+6q^2+20q^3+50q^4+104q^5+190q^6+313q^7+473q^8+663q^9+868q^{10}\\[8pt]
&+1068q^{11}+1240q^{12}+1362q^{13}+1417q^{14}+1398q^{15}+1307q^{16}\\[8pt]
&+1157q^{17}+968q^{18}+763q^{19}+564q^{20}+388q^{21}+247q^{22}\\[8pt]
&+143q^{23}+74q^{24}+33q^{25}+12q^{26}+3q^{27}+q^{28},\\[8pt]
d_9(q)=q&+7q^2+27q^3+77q^4+181q^5+371q^6+684q^7+1157q^8+1820q^9\\[8pt]
&+2687q^{10}+3749q^{11}+4969q^{12}+6281q^{13}+7594q^{14}+8802q^{15}\\[8pt]
&+9796q^{16}+10480q^{17}+10785q^{18}+10680q^{19}+10176q^{20}+9324q^{21}\\[8pt]
&+8209q^{22}+6935q^{23}+5611q^{24}+4337q^{25}+3192q^{26}+
2227q^{27}+1465q^{28}\\[8pt]
&+901q^{29}+513q^{30}+266q^{31}+123q^{32}+49q^{33}+16q^{34}+4q^{35}.
 \end{align*}
Observe that  $d_9(q)$ has the following spiral property:
\[ 1 < {\underline 4} < 7 < \underline{16}  < 27 < \underline{49}
 < 77 < \underline{123}  < 181 < \underline{266}  < 371
  <  \underline{513}
 < 684 < \underline{901} < \cdots
.\]
 The underlined numbers are chosen backwards
starting from the last coefficient. In general, a spiral sequence
consists of two interlacing sequences. Moreover, the spiral property
implies unimodality. Our objective is to find two monotone ratio
sequences that imply the spiral property and the log-concavity. For
$d_9(q)$, the two monotone sequences are given below:
\[  {1 \over 4} < {7 \over 16} < {27 \over 49} <  {77 \over 123}
 <  {181\over 266 } < {371\over 513} < {684\over 901} < \cdots
 <{10480\over
 10680}< 1, \]
\[  {4\over 7} < {16\over  27}  < {49\over 77} <
{123\over 181} < {266\over 371}
  <  {513\over 684} < {901\over 1157}< \cdots <{10680\over 10785} < 1
.\] For $d_8(q)$, since the first term and the last term are both
equal to one, we need to ignore the last term in order to construct
 two monotone ratio sequences which are given below:
\[  {1 \over 3} < {6 \over 12} < {20 \over 33} <  {50 \over 74}
 <  {104\over 143 } < {190\over 247} < {313\over 388} < \cdots <
 {1362\over 1398}< 1, \]
\[  {3\over 6} < {12\over  20}  < {33\over 50} <
{74\over 104} < {143\over 190}
  <  {247\over 313} < {388\over 473}< \cdots <{1398\over 1417}< 1
.\] It is easily seen that the spiral property can be recovered from
the above ratio monotone property.

To conclude the introduction, we remark that the above ratio
monotone property implies the log-concavity. As shown in the above
examples, for $n=8$ the $q$-derangement number is log-concave except
for the last term, whereas $d_9(q)$ is log-concave in the usual
sense. In fact, this fact is valid in general for $n$ being even and
odd.

We say that a positive sequence $a_1, a_2, \ldots, a_n$ satisfies
the ratio monotone property if
\begin{align*}
\frac{a_1}{a_n}&<\frac{a_2}{a_{n-1}}<
\cdots<\frac{a_{i-1}}{a_{n+2-i}}
<\frac{a_i}{a_{n+1-i}}<\cdots,\\[8pt]
\frac{a_n}{a_2}&<\frac{a_{n-1}}{a_3}<\cdots <\frac{a_{n+2-i}}{a_i}<
\frac{a_{n+1-i}}{a_{i+1}}<\cdots.
\end{align*}
For comparison, we recall that a positive sequence $a_1, a_2,
\ldots, a_n$  is log-concave if
\[
\frac{a_1}{a_2}<\frac{a_2}{a_3}< \cdots <\frac{a_i}{a_{i+1}}<
\cdots.
\]
Note that the ratio monotone property implies the log-concavity.
This is because the two inequalities
 \[
\frac{a_{i-1}}{a_{n+2-i}}<\frac{a_i}{a_{n+1-i}},\ \ \ \ \
\frac{a_{n+2-i}}{a_i}<\frac{a_{n+1-i}}{a_{i+1}},
 \]
when written as
\[
a_{i-1}a_{n+1-i}<a_{n+2-i}a_i, \ \ \ a_{n+2-i}a_{i+1}<a_ia_{n+1-i},
\]
yields the relation
$a_{i-1}a_{n+1-i}a_{n+2-i}a_{i+1}<a_{n+2-i}a_ia_ia_{n+1-i}$, which
implies
\[
\frac{a_{i-1}}{a_i}<\frac{a_i}{a_{i+1}}.
\]
Thus the ratio monotone property of $d_n(q)$ leads to the
aforementioned log-concavity.

\section{The Main Theorem}

Let $\beta_n$ denote the degree of $d_n(q)$, $n\geq 2$. It is easily
seen that
\begin{align*}
\beta_n & = \left\{
     \begin{array}{ll}
         {n \choose 2}, \qquad & \mbox{if $n$ is even}, \\ [8pt]
        {n \choose 2}-1, \qquad
              & \mbox{if $n$ is odd}.
     \end{array} \right.
\end{align*}

Set
\[ d_n(q)=A_n(1)q+A_n(2)q^2+\cdots+A_n(\beta_n)q^{\beta_n}.\]
The ratio monotone property for $d_n(q)$ $n\geq 6$ can be stated in
the following theorem. It turns out the structure of the ratio
sequences depends on the residue of $n$ modulo $4$.

\begin{thm}{\label{theorem}}
For $n \geq 6$, we have the following monotone ratio sequences
formed by the coefficients of $d_n(q)$:
\begin{itemize}
\item[1.]
$n\equiv 1\  {\rm mod}  \ 4$: Let $r={n(n-1)\over 4}$. Then we have
\begin{align}\label{thin1}
\frac{A_n(1)}{A_n(\beta_n)}< \frac{A_n(2)}{A_n(\beta_{n}-1)}<\cdots
<\frac{A_n(r-1)}{A_n(r+1)}<1
\end{align}
and
\begin{align}\label{thin2}
\frac{A_n(\beta_n)}{A_n(2)}< \frac{A_n(\beta_n-1)}{A_n(3)}<\cdots <
\frac{A_n(r+1)}{A_n(r)}<1.
\end{align}
\item[2.]
If $n\equiv 2 \ {\rm mod}  \ 4$: Let $r={n(n-1)-2\over 4}$. Then we
have
\begin{align}\label{thin3}
\frac{A_n(1)}{A_n(\beta_n-1)}< \frac{A_n(2)}{A_n(\beta_n-2)}< \cdots
<\frac{A_n(r)}{A_n(r+1)}<1
\end{align}
and
\begin{align}\label{thin4}
\frac{A_n(\beta_n-1)}{A_n(2)}< \frac{A_n(\beta_n-2)}{A_n(3)}< \cdots
<\frac{A_n(r+2)}{A_n(r)}< 1.
\end{align}
\item[3.]
If $n\equiv 3 \ {\rm mod}  \ 4$: Let $r={n(n-1)-2\over 4}$. Then we
have
\begin{align}\label{thin5}
\frac{A_n(1)}{A_n(\beta_n)}< \frac{A_n(2)}{A_n(\beta_n-1)}< \cdots
<\frac{A_n(r)}{A_n(r+1)}<1
\end{align}
and
\begin{align}\label{thin6}
\frac{A_n(\beta_n)}{A_n(2)}< \frac{A_n(\beta_n-1)}{A_n(3)}< \cdots
<\frac{A_n(r+2)}{A_n(r)}<1.
\end{align}
\item[4.]
If $n\equiv 0 \ {\rm mod}  \ 4$: Let $r={n(n-1)\over 4}$. Then we
have
\begin{align}\label{thin7}
\frac{A_n(1)}{A_n(\beta_n-1)}< \frac{A_n(2)}{A_n(\beta_n-2)}< \cdots
<\frac{A_n(r-1)}{A_n(r+1)}< 1
\end{align}
and
\begin{align}\label{thin8}
\frac{A_n(\beta_n-1)}{A_n(2)}< \frac{A_n(\beta_n-2)}{A_n(3)}< \cdots
<\frac{A_n(r+1)}{A_n(r)}<1 .
\end{align}
\end{itemize}
\end{thm}

To prove Theorem \ref{theorem}, we need the following two lemmas.

\begin{lem}\label{lemma1}
Suppose that $a_1, a_2, \ldots, a_{k+1}, b_1, b_2, \ldots, b_{k+1}$
 $(k\geq 1)$ are positive real numbers satisfying
\begin{align}\label{assum}
\frac{b_1}{a_1}<\frac{b_2}{a_2}<\cdots
<\frac{b_k}{a_k}<\frac{b_{k+1}}{a_{k+1}}.
\end{align}
Then the following inequalities hold:
\begin{align}
\frac{b_1+b_2+\cdots +b_k}{a_1+a_2+\cdots +a_k+1}&<
\frac{b_1+b_2+\cdots +b_k+b_{k+1}}{a_1+a_2+\cdots +a_k+a_{k+1}+1},
\label{ineq1}\\[8pt]
\frac{b_1+b_2+\cdots +b_k}{a_1+a_2+\cdots +a_k}&<
\frac{b_1+b_2+\cdots +b_k+b_{k+1}}{a_1+a_2+\cdots +a_k+a_{k+1}},
\label{ineq2} \\[8pt]
 \frac{b_1+b_2+\cdots +b_k}{a_1+a_2+\cdots +a_k}& <
\frac{b_2+\cdots +b_k+b_{k+1}}{a_2+\cdots +a_k+a_{k+1}}.
\label{ineq3}
\end{align}
\end{lem}

\noindent
 {\it Proof.} From  the conditions \eqref{assum}, we have
\[
a_ib_{k+1}>b_ia_{k+1},\ \ a_1b_{i+1}>b_1a_{i+1}, \qquad 1\leq i\leq
k.
\]
Consequently,
\begin{align}
(a_1+&\cdots +a_k)(b_1+\cdots
+b_k)+(b_1+\cdots+b_k)+b_{k+1}+b_{k+1}(a_1+\cdots+a_k)\nonumber\\[8pt]
&>(a_1+ \cdots +a_k)(b_1+\cdots +b_k)+a_{k+1}(b_1+\cdots
+b_k)+(b_1+\cdots +b_k),
 \label{ineq4}\\[8pt]
 (b_1+&b_2+\cdots+b_k)(a_1+a_2+\cdots
 +a_k)+(a_1+a_2+\cdots+a_k)b_{k+1} \nonumber\\[8pt]
 &>(b_1+b_2+\cdots+b_k)(a_1+a_2+\cdots
 +a_k)+(b_1+b_2+\cdots+b_k)a_{k+1} \label{ineq5}
 \end{align}
 and
 \begin{align}
a_1(b_2+&\cdots +b_{k+1})+(a_2+\cdots+a_k)(b_2+\cdots +b_k)
+b_{k+1}(a_2+\cdots+a_k)\nonumber \\
&>b_1(a_2+\cdots +a_{k+1})+(a_2+\cdots+a_k)(b_2+\cdots
+b_k)+a_{k+1}(b_2+\cdots +b_k). \label{ineq6}
\end{align}
Therefore,  \eqref{ineq1}, \eqref{ineq2} and \eqref{ineq3} can be
derived from  \eqref{ineq4}, \eqref{ineq5} and \eqref{ineq6},
respectively. \qed

\begin{lem}\label{lemma2}
If $m$ is even and $m \geq 6$, then
\begin{align*}
A_m(\beta_m-1)&=\frac{m}{2}-1,\ A_m(2)=m-2,\
A_m(3)=\frac{m(m-3)}{2},\
A_m(\beta_m-2)=\frac{m^2}{4}-\frac{m}{2},\\[8pt]
A_m(4)&=\frac{(m-2)(m+2)(m-3)}{6},\ \ \
A_m(\beta_m-3)=\frac{1}{12}m^3-\frac{1}{8}m^2-\frac{1}{12}m-1.
\end{align*}
\end{lem}

\noindent {\it Proof.} We use induction on $m$. It is evident that
the lemma holds for $m=6$. Moreover, $A_j(1)=1$ for $j \geq 2$.
Suppose that the lemma is true for $t \geq 6$, where $t$ is even. We
consider the case $m= t +2$. The formula \eqref{recursion} on
$d_n(q)$ leads to  the recurrence relation for $n\geq 1$:
\begin{eqnarray}\label{re}
A_{n+1}(k)=\left\{
     \begin{array}{lllll}
    \sum\limits_{i=1}^k A_n(i), & \mbox{ $1\leq k \leq n+1$ },
     \\[12pt]
       \sum\limits_{i=k-n}^k A_n(i),
              & \mbox{ $n+1\leq k\leq \beta_n$ },\\[12pt]
    \sum\limits_{i=k-n}^{\beta_n} A_n(i),
              & \mbox{ $\beta_n\leq k < \beta_{n+1}$ },\\[12pt]
              A_n(\beta_n-1) +A_n(\beta_n),
              &\mbox{ $k =\beta_{n+1}$ and $n$ is even},\\[12pt]
              1,
              &\mbox{ $k =\beta_{n+1}$ and $n$ is odd}.
     \end{array} \right.
\end{eqnarray}

Using the above recursion, by induction,  we have
\begin{align*}
A_{t+2}(2)=A_{t+1}(2)+A_{t+1}(1)=A_{t}(2)
+A_{t}(1)+A_{t+1}(1) =t-2+1+1=t,\\[8pt]
A_{t+2}(\beta_{t+2}-1)=A_{t+1}(\beta_{t +1})
=A_{t}(\beta_{t})+A_{t}(\beta_{t}-1)=1+\frac{t}{2}
-1=\frac{t+2}{2}-1
\end{align*}
and
\begin{align*}
A_{t+2}(3)&=A_{t+2}(2)+A_{t+1}(3)=A_{t+2}(2)+
A_{t}(3)+A_{t}(2)+A_{t}(1)\\[8pt]
&=\frac{(t+2)(t+2-3)}{2}.
\end{align*}
Moreover,
\begin{align*}
A_{t+2}(\beta_{t+2}-2)&=A_{t+1}(\beta_{t+1}-1)+
A_{t+1}(\beta_{t+1})\\[8pt]
&=A_{t}(\beta_{t}-2)+ A_{t}(\beta_{t}-1)+A_{t}(\beta_{t})+
A_{t}(\beta_{t}-1)+A_{t}(\beta_{t})\\[8pt]
&=\frac{(t+2)^2}{4}-\frac{(t+2)}{2}.
\end{align*}
Analogously, we can compute $A_{t+2}(4)$ and
$A_{t+2}(\beta_{t+2}-3)$.
 This
completes the proof. \qed

We are now ready to present the proof of the main theorem.

\noindent
 {\it Proof of Theorem \ref{theorem}.}
 We use induction on $n$. For $n=6,7,8,9$, it is easy to verify that the
theorem is true. Here we consider only two cases: $m\equiv 0 \mod 4$
and $m\equiv 1 \mod 4$, namely, only the sequences \eqref{thin1} and
\eqref{thin2}, \eqref{thin3} and \eqref{thin4}.
 The proofs of \eqref{thin5} and
\eqref{thin6} are the similar to those for  \eqref{thin1} and
\eqref{thin2}, and  the proofs of \eqref{thin7} and \eqref{thin8}
are the similar to those for \eqref{thin3} and \eqref{thin4}.

 Suppose that the theorem holds for
$m$, where $m\equiv 0\ {\rm mod}  \ 4$, namely,
\begin{align}\label{induction1}
\frac{A_m(1)}{A_m(\beta_m-1)}&<\frac{A_m(2)}{A_m(\beta_m-2)}< \cdots
 <\frac{A_m(\frac{m(m-1)}{4}-1)}{A_m(\frac{m(m-1)}{4}+1)}<1.
\end{align}
and
\begin{align}\label{indu2}
\frac{A_m(\beta_m-1)}{A_m(2)}&< \frac{A_m(\beta_m-2)}{A_m(3)}<
\cdots <\frac{A_m(\frac{m(m-1)}{4}+1)}{A_m(\frac{m(m-1)}{4})}<1.
\end{align}
 We now proceed to show that the theorem also holds for
 $n=m+1$. Let $r={m(m+1) \over 4}$. The desired monotone ratio
 sequences are stated as follows:
\begin{align}\label{thin1-m+1}
\frac{A_{m+1}(1)}{A_{m+1}(\beta_{m+1})}<
\frac{A_{m+1}(2)}{A_{m+1}(\beta_{m+1}-1)}<\cdots
<\frac{A_{m+1}(r-1)}{A_{m+1}(r+1)}< 1
\end{align}
and
\begin{align}\label{thin2-m+1}
\frac{A_{m+1}(\beta_{m+1})}{A_{m+1}(2)}<
\frac{A_{m+1}(\beta_{m+1}-1)}{A_{m+1}(3)}<\cdots <
\frac{A_{m+1}(r+1)}{A_{m+1}(r)}<1.
\end{align}

Since $m\equiv 0\ {\rm mod}  \ 4$, we have $\beta_m= {m(m-1)\over
2}$ and $\beta_{m+1} = {m(m+1)\over 2} -1$. We now aim to prove
\eqref{thin1-m+1}. We will divide the ratio sequence
  \eqref{thin1-m+1} into three segments. First,  for $1 \leq k\leq
m-1$, from \eqref{re},  we have
\begin{align}
A_{m+1}(k)&= \sum_{i=1}^k A_m(i), \label{case1}\\[8pt]
A_{m+1}(k+1)&=\sum_{i=1}^{k+1}A_m(i). \label{case1-1}
 \end{align}
It is easily checked that $ \beta_m+1
 \leq \beta_{m+1}+1-k \leq   \beta_{m+1}$. Note that $A_m(
 \beta_m)=1$, and we claim that
\begin{align}
A_{m+1}(\beta_{m+1}+1-k)&= \sum_{i=1}^kA_m(\beta_m-i)+1, \label{case1-2}\\[8pt]
 A_{m+1}(\beta_{m+1}-k)&=
\sum_{i=1}^{k+1}A_m(\beta_m-i)+1.\label{case1-3}
\end{align}
Clearly, for $1<k \leq m-1$, (\ref{case1-2}) is given by the above
recurrence relation \eqref{re}.  It suffices to check that
(\ref{case1-2}) holds for $k=1$. In this case,
\[ A_{m+1}(\beta_{m+1})=A_m(\beta_m-1) +
A_m(\beta_m)=A_m(\beta_m-1)+1,\]
 which is in agreement with the
$k=1$ case of (\ref{case1-2}).

Because of \eqref{ineq1} and \eqref{induction1},
 the following relation holds
\begin{align*}
\frac
{\sum\limits_{i=1}^kA_m(i)}{\sum\limits_{i=1}^kA_m(\beta_m-i)+1}
 < \frac{\sum\limits_{i=1}^{k+1}A_m(i)}
{\sum\limits_{i=1}^{k+1}A_m(\beta_m-i)+1},
\end{align*}
which can be recast as
\begin{align}\label{k1}
\frac{A_{m+1}(k)}{A_{m+1}(\beta_{m+1}+1-k)}<\frac{A_{m+1}(k+1)}
{A_{m+1}(\beta_{m+1}-k)},\
 \ \ 1\leq k  \leq m-1.
\end{align}

Second, for $k=m$, from the recurrence relation \eqref{re},  we get
\begin{align*}
A_{m+1}(m)&=\sum_{i=1}^mA_m(i),\\[8pt]
A_{m+1}(m+1)&=\sum_{i=1}^{m+1}A_m(i),\\[8pt]
A_{m+1}(\beta_{m+1}+1-m)&=1+\sum_{i=1}^mA_m(\beta_m-i),\\[8pt]
A_{m+1}(\beta_{m+1}-m)&=\sum_{i=1}^{m+1}A_m(\beta_m-i).
\end{align*}
From \eqref{ineq2} and \eqref{induction1} it can be deduced that
\begin{align*}
\frac{\sum\limits_{i=1}^mA_m(i)}
{1+\sum\limits_{i=1}^mA_m(\beta_m-i)}
<\frac{\sum\limits_{i=1}^mA_m(i)}
{\sum\limits_{i=1}^mA_m(\beta_m-i)}
<\frac{\sum\limits_{i=1}^{m+1}A_m(i)}
{\sum\limits_{i=1}^{m+1}A_m(\beta_m-i)},
\end{align*}
which can be restated as
\begin{align}\label{k2}
\frac{A_{m+1}(m)}{A_{m+1}(\beta_{m+1}+1-m)}<
\frac{A_{m+1}(m+1)}{A_{m+1}(\beta_{m+1}-m)}.
\end{align}
Finally, for $m<k \leq r-1$,  we have $r+1\leq \beta_{m+1}+1-k \leq
\beta_m-1$,
   the recurrence relation \eqref{re}
implies that
 \begin{align}
A_{m+1}(k)&=\sum_{i=0}^mA_m(k-i),\label{case2-1}\\[8pt]
A_{m+1}(k+1)&=\sum_{i=0}^mA_m(k+1-i) \label{case2-2}
\end{align}
and
 \begin{align}
A_{m+1}(\beta_{m+1}-k)&=\sum_{i=0}^mA_m(\beta_m-k-1+i), \label{case2-3}\\[8pt]
A_{m+1}(\beta_{m+1}+1-k)&=\sum_{i=0}^mA_m(\beta_m-k+i).
\label{case2-4}
\end{align}
We have from  \eqref{ineq3}  and  \eqref{induction1},
\begin{align}\label{k3}
\frac{\sum\limits_{k=0}^mA_m(k-i)}
{\sum\limits_{i=0}^mA_m(\beta_m-k+i)}
<\frac{\sum\limits_{k=0}^mA_m(k+1-i)}
{\sum\limits_{i=0}^mA_m(\beta_m-k-1+i)} .
\end{align}
Taking into account altogether \eqref{k1}, \eqref{k2} and
\eqref{k3}, we conclude that  for $1\leq k \leq r-1$,
\begin{equation} \label{amb}
\frac{A_{m+1}(k)}{A_{m+1}(\beta_{m+1}+1-k)}< \frac
{A_{m+1}(k+1)}{A_{m+1}(\beta_{m+1}-k)}.
\end{equation}
Plugging  $k =r-1$ in  \eqref{amb} gives
\[\frac{A_{m+1}(r-1)}{A_{m+1}(r+1)}<\frac{A_{m+1}(r)}{A_{m+1}(r)}
=1.
\]
Thus we have established \eqref{thin1-m+1}.

 Now, we turn to the ratio sequence \eqref{thin2-m+1}. We
   also need to consider three cases. First, when $2 \leq k <m+1$,
  then $\beta_m \leq \beta_{m+1}+1-k \leq \beta_{m+1}-1$, and
  $A_{m+1}(k)$, $A_{m+1}(k+1)$, $A_{m+1}(\beta_{m+1}+1-k)$ can be
  expressed in terms of $A_m(i)$
   by \eqref{case1}, \eqref{case1-1}, \eqref{case1-2},
   respectively. Combining the two cases $k=2$ and $2<k<m+1$, we
   reach the assertion
\[
A_{m+1}(\beta_{m+1}+2-k)=1+\sum_{i=1}^{k-1}A_m(\beta_m-i).
\]
On the other hand, by Lemma \ref{lemma2}, we find that
\begin{align*}
A_m&(\beta_m-2)\left(1+A_m(2)\right)
-A_m(3)\left(1+A_m(\beta_m-1)\right)\\[8pt]
&=\left(\frac{m^2}{4}-\frac{m}{2}\right)(1+m-2)-\frac{m(m-3)}{2}
\left(1+\frac{m}{2}-1\right) =\frac{m}{2}>0.
\end{align*}
Furthermore,  the inductive hypothesis \eqref{indu2} implies that
\[\frac{A_m(\beta_m-k)}{A_m(k+1)}\geq \frac{A_m(\beta_m-2)}{A_m(3)}
>\frac{1+A_m(\beta_m-1)}{1+A_m(2)}.
\]
Thus,
\begin{align}\label{l1}
A_m(k+1)\left(1+A_m(\beta_m-1)\right)<A_m(\beta_m-k)\left(1+A_m(2)\right).
\end{align}
Again, as a consequence of \eqref{indu2},
 we see that
\begin{align}\label{l2}
 A_m(k+1)A_m(\beta_m+1-i)<A_m(\beta_m-k)A_m(i),  \ \ 3 \leq i \leq k.
\end{align}
In light of \eqref{l1} and \eqref{l2}, we obtain the following
relation
\begin{align*}
&\left(1+\sum\limits_{i=2}^kA_m(i)\right)
\left(1+\sum\limits_{i=1}^{k-1}A_m(\beta_m-i)\right)\\[8pt]
& \qquad +A_m(k+1)\left(1+A_m(\beta_m-1)\right) +A_m(k+1)
\left(\sum\limits_{i=2}^{k-1}A_m(\beta_m-i)\right)\\[8pt]
& \quad <\left(1+\sum\limits_{i=2}^kA_m(i)\right)
\left(1+\sum\limits_{i=1}^{k-1}A_m(\beta_m-i)\right)\\[8pt]
&  \qquad +A_m(\beta_m-k)\left(1+A_m(2)\right)+
A_m(\beta_m-k)\left(\sum\limits_{i=3}^kA_m(i)\right).
\end{align*}
Therefore,
\begin{align*}
\frac{1+\sum\limits_{i=1}^{k-1}A_m(\beta_m-i)}
{1+\sum\limits_{i=2}^kA_m(i)}
<\frac{1+\sum\limits_{i=1}^kA_m(\beta_m-i)}
{1+\sum\limits_{i=2}^{k+1}A_m(i)},
\end{align*}
which can be restated as
\begin{align}\label{k4}
\frac{A_{m+1}(\beta_{m+1}+2-k)}{A_{m+1}(k)}< \frac
{A_{m+1}(\beta_{m+1}+1-k)}{A_{m+1}(k+1)},\ \ \ \ 2\leq k < m+1.
\end{align}
Second, we consider the case $k=m+1$. We have from
\eqref{induction1},
 \begin{align}\label{leq1}
 A_m(i)<A_{m}(\beta_m-i), \ \ \ \ 2 \leq i \leq m+1.
 \end{align}
Note that the following relation is implied by \eqref{indu2},
\begin{align}\label{leq2}
 A_m(m+2)A_m(\beta_m-i)<A_{m}(\beta_m-m-1)A_{m}(i+1), \ \ \ \ 3 \leq i \leq m.
 \end{align}
In addition, Lemma \ref{lemma2} enables us to check
\begin{align*}
A_m(\beta_m-3)&\left(A_m(2)+A_m(3)\right)-A_m(4)
\left(1+A_m(\beta_m-1)+A_m(\beta_m-2)\right)\\[8pt]
=&\left(\frac{m^3}{12}-\frac{m^2}{8}-\frac{m}{12}-1\right)
\left(m-2+\frac{m(m-3)}{2}\right)\\[8pt]
& -\frac{(m-2)(m+2)(m-3)}{6} \left(
1+\frac{m}{2}-1+\frac{m^2}{4}-\frac{m}{2}\right)\\[8pt]
=&\frac{(m-4)(m^3+5m^2-14m-24)}{48},
\end{align*}
which is positive for $m \geq 6$. This implies that
\[
\frac{A_m(\beta_m-3)}{A_m(4)}>
\frac{1+A_m(\beta_m-1)+A_m(\beta_m-2)}{A_m(2)+A_m(3)}.
\]
Note that the sequence \eqref{indu2} contains the following
inequality
\[
\frac{A_m(\beta_m-m-1)}{A_{m}(m+2)}>\frac{A_m(\beta_m-3)}{A_m(4)}.\]
Hence
\[
\frac{A_m(\beta_m-m-1)}{A_m(m+2)}>
\frac{1+A_m(\beta_m-1)+A_m(\beta_m-2)}{A_m(2)+A_m(3)},
\]
or, equivalently,
\begin{align}\label{leq3}
&\left(A_m(2)+A_m(3)\right)A_m(\beta_m-m-1) \nonumber\\[8pt]
&\ \ \ \ \ \ \ \ >A_m(m+2)
\left(1+A_m(\beta_m-1)+A_m(\beta_m-2)\right).
\end{align}
From \eqref{leq1}, \eqref{leq2} and \eqref{leq3}, we deduce that
\begin{align*}
&\left(\sum\limits_{i=1}^mA_m(\beta_m-i)\right) \left(
\sum\limits_{i=2}^{m+1}A_m(i) \right)
+A_m(\beta_m-m-1)\left(A_m(2)+A_m(3)\right)\\[8pt]
&\qquad +A_m(\beta_m-m-1) \left(\sum\limits_{i=4}^{m+1}A_m(i)\right)
+\sum\limits_{i=1}^{m+1}A_m(\beta_m-i)\\[8pt]
&\quad
>\left(\sum\limits_{i=1}^mA_m(\beta_m-i)\right)\left(
\sum\limits_{i=2}^{m+1}A_m(i) \right)
+A_m(m+2)\left(1+A_m(\beta_m-1)+A_m(\beta_m-2)\right)
\\[8pt]
 &\qquad + A_m(m+2) \left( \sum\limits_{i=3}^{m}
 A_m(\beta_m-i) \right) +
\sum\limits_{i=2}^{m+1}A_m(i).
\end{align*}
Thus,
\begin{align*}
\frac{\sum\limits_{i=1}^mA_m(\beta_m-i)+1}{\sum\limits_{i=2}^{m+1}
A_m(i)+1} <\frac{\sum\limits_{i=1}^{m+1}A_m(\beta_m-i)}
{\sum\limits_{i=2}^{m+2}A_m(i)},
\end{align*}
which can be expressed in terms of $A_{m+1}(k)$ as follows:
\begin{align}\label{k5}
\frac{A_{m+1}(\beta_{m+1}+1-m)}{A_{m+1}(m+1)}< \frac
{A_{m+1}(\beta_{m+1}-m)}{A_{m+1}(m+2)}.
\end{align}

 Next, we come to  the range $m+1 < k\leq r$. In this case, we
 have
  $\frac{m(m+1)}{4}+1
 \leq \beta_{m+1}+2-k \leq \beta_m -1$. Hence $A_{m+1}(k)$,
  $A_{m+1}(k+1)$, $A_{m+1}(\beta_{m+1}+1-k)$ satisfy the
  recurrence relations
   \eqref{case2-1}, \eqref{case2-2}, \eqref{case2-4}, respectively.
   Moreover,
\[
A_{m+1}(\beta_{m+1}+2-k)=\sum_{i=0}^mA_{m}(\beta_m+1-k+i).
\]
Again,  justified by \eqref{ineq3} and \eqref{indu2},
  we obtain
\begin{align}\label{k6}
\frac{\sum\limits_{i=0}^mA_m(\beta_m+1-k+i)}{\sum\limits_{i=0}^m
A_m(k-i)} <\frac{\sum\limits_{i=0}^mA_m(\beta_m-k+i)}
{\sum\limits_{i=0}^mA_m(k+1-i)} .
\end{align}
Combining \eqref{k4}, \eqref{k5} and \eqref{k6}, we find
\[
\frac{A_{m+1}(\beta_{m+1}+2-k)}{A_{m+1}(k)}< \frac
{A_{m+1}(\beta_{m+1}+1-k)}{A_{m+1}(k+1)},\ \ \ \ 2\leq k \leq r.
\]
Setting $k=r$ in the above inequality gives
\[
\frac{A_{m+1}(r+1)}{A_{m+1}(r)}<\frac{A_{m+1}(r)}{A_{m+1}(r+1)},
\]
which immediately becomes
\[
\frac{A_{m+1}(r+1)}{A_{m+1}(r)}<1.
\]
Thus the proof of  \eqref{thin2-m+1} is complete.

We now attempt to prove  \eqref{thin3} and \eqref{thin4}. Suppose
that the theorem is valid for $m$, where $m\equiv 1 \ {\rm mod} \
4$, i.e.,
\begin{align}\label{indu3}
\frac{A_m(1)}{A_m(\beta_m)}< \frac{A_m(2)}{A_m(\beta_{m}-1)}<\cdots
<\frac{A_m(\frac{m(m-1)}{4}-1)}{A_m(\frac{m(m-1)}{4}+1)}<1
\end{align}
and
\begin{align}\label{indu4}
\frac{A_m(\beta_m)}{A_m(2)}< \frac{A_m(\beta_m-1)}{A_m(3)}<\cdots <
\frac{A_m(\frac{m(m-1)}{4}+1)}{A_m(\frac{m(m-1)}{4})}<1.
\end{align}
 We will show that the theorem also holds for
 $n=m+1$. Let $r=\frac{m(m+1)-2}{4}$. Note that in this case
  $\beta_m={m\choose 2}-1$ and $\beta_{m+1}={m+1\choose 2}$.
  Our goal is to prove  the
 following two relations:
\begin{align}\label{want1}
\frac{A_{m+1}(1)}{A_{m+1}(\beta_{m+1}-1)}<
\frac{A_{m+1}(2)}{A_{m+1}(\beta_{m+1}-2)}<\cdots
<\frac{A_{m+1}(r)}{A_{m+1}(r+1)}<1
\end{align}
and
\begin{align}\label{want2}
\frac{A_{m+1}(\beta_{m+1}-1)}{A_{m+1}(2)}<
\frac{A_{m+1}(\beta_{m+1}-2)}{A_{m+1}(3)}< \cdots
<\frac{A_{m+1}(r+2)}{A_{m+1}(r)}< 1.
\end{align}

We first consider \eqref{want1}. When $1\leq k \leq m$,
 then $\beta_m+1\leq \beta_{m+1}-k \leq \beta_{m+1}-1$,
 $A_{m+1}(k)$ and $A_{m+1}(k+1)$ can be expressed by \eqref{case1}
and \eqref{case1-1}, respectively. From \eqref{re} it follows that
\begin{align*}
A_{m+1}(\beta_{m+1}-k)&=\sum_{i=1}^{k}A_m(\beta_m+1-i),\\[8pt]
A_{m+1}(\beta_{m+1}-k-1)&=\sum_{i=1}^{k+1}A_m(\beta_m+1-i).
\end{align*}
Using  \eqref{indu3} and \eqref{ineq2}, we obtain
\begin{align*}
\frac{\sum\limits_{i=1}^kA_m(i)}
{\sum\limits_{i=1}^{k}A_m(\beta_m+1-i)} <
\frac{\sum\limits_{i=1}^{k+1}A_m(i)}
{\sum\limits_{i=1}^{k+1}A_m(\beta_m+1-i)},
\end{align*}
which can be rewritten as
\begin{align}\label{w-1}
\frac{A_{m+1}(k)}{A_{m+1}(\beta_{m+1}-k)}
<\frac{A_{m+1}(k+1)}{A_{m+1}(\beta_{m+1}-k-1)},\ \ \ \ 1 \leq k \leq
m.
\end{align}
When $m<k \leq r$, then $r+1 \leq \beta_{m+1}-k \leq \beta_m$,
$A_{m+1}(k)$ and $A_{m+1}(k+1)$
 can be expressed by \eqref{case2-1} and \eqref{case2-2},
 respectively. From \eqref{re}, we have
 \begin{align*}
A_{m+1}(\beta_{m+1}-k)&=\sum_{i=0}^mA_m(\beta_m+1-k+i),\\[8pt]
A_{m+1}(\beta_{m+1}-k-1)&=\sum_{i=0}^mA_m(\beta_m-k+i).
 \end{align*}
Now, based on \eqref{indu3} and \eqref{ineq3}, we derive
 \begin{align}\label{w-2}
\frac{\sum\limits_{i=0}^mA_m(k-i)}
{\sum\limits_{i=0}^mA_m(\beta_m+1-k+i)}<
\frac{\sum\limits_{i=0}^mA_m(k+1-i)}
{\sum\limits_{i=0}^mA_m(\beta_m-k+i)}.
 \end{align}
Combing \eqref{w-1} and \eqref{w-2} leads to
\begin{align}\label{w-3}
\frac{A_{m+1}(k)}{A_{m+1}(\beta_{m+1}-k)}
<\frac{A_{m+1}(k+1)}{A_{m+1}(\beta_{m+1}-k-1)},\ \ \ \ 1 \leq k \leq
r.
\end{align}
Substituting $k=r$ into \eqref{w-3}, we get
\begin{align*}
\frac{A_{m+1}(r)}{A_{m+1}(r+1)} <\frac{A_{m+1}(r+1)}{A_{m+1}(r)},
\end{align*}
which implies
\[
\frac{A_{m+1}(r)}{A_{m+1}(r+1)}<1.
\]
So we have shown that \eqref{want1} is valid.

It remains to prove \eqref{want2}. We still need to consider three
cases for the index $k$. First, when $ 2 \leq k \leq m$, then
$\beta_m+1 \leq \beta_{m+1}-k \leq \beta_{m+1}-2$,
 $A_{m+1}(k)$ and $A_{m+1}(k+1)$ can be expressed by \eqref{case1}
and \eqref{case1-1}, respectively. In view of \eqref{re}, we see
that
\begin{align*}
A_{m+1}(\beta_{m+1}+1-k)&=\sum_{i=2}^{k}A_m(\beta_m+2-i),\\[8pt]
A_{m+1}(\beta_{m+1}-k)&=\sum_{i=2}^{k+1}A_m(\beta_m+2-i).
\end{align*}
Note that $A_m(1)=1$. Using \eqref{ineq1} and \eqref{indu4}, we
deduce that
\begin{align*}
\frac{\sum\limits_{i=2}^{k}A_m(\beta_m+2-i)}
{1+\sum\limits_{i=2}^kA_m(i)} <
\frac{\sum\limits_{i=2}^{k+1}A_m(\beta_m+2-i)}
{1+\sum\limits_{i=2}^{k+1}A_m(i)},
\end{align*}
which  yields
\begin{align}\label{g1}
\frac{A_{m+1}(\beta_{m+1}+1-k)}{A_{m+1}(k)}<
\frac{A_{m+1}(\beta_{m+1}-k)}{A_{m+1}(k+1)}, \ \ \ \ 2 \leq k \leq
m.
\end{align}
Second, when $k=m+1$, from \eqref{re} we get
\begin{align*}
A_{m+1}(\beta_{m+1}-m-1)&=\sum_{i=2}^{m+2}A_m(\beta_m+2-i),\\[8pt]
A_{m+1}(m+2)&=\sum_{i=2}^{m+2}A_m(i).
\end{align*}
Using  \eqref{indu4} and \eqref{ineq2}, we find
\begin{align*}
\frac{\sum\limits_{i=2}^{m+1}A_m(\beta_m+2-i)}
{1+\sum\limits_{i=2}^{m+1}A_m(i)}  <
\frac{\sum\limits_{i=2}^{m+1}A_m(\beta_m+2-i)}
 {\sum\limits_{i=2}^{m+1}A_m(i)}
< \frac{\sum\limits_{i=2}^{m+2}A_m(\beta_m+2-i)}
{\sum\limits_{i=2}^{m+2}A_m(i)},
\end{align*}
which can be restated as
\begin{align}\label{g2}
\frac{A_{m+1}(\beta_{m+1}-m)}{A_{m+1}(m+1)}<
\frac{A_{m+1}(\beta_{m+1}-m-1)}{A_{m+1}(m+2)}.
\end{align}
Finally, when $m+1< k \leq r$,  $r+1 \leq \beta_{m+1}-k \leq
\beta_m-1$, $A_{m+1}(k)$ and $A_{m+1}(k+1)$
 can be expressed by \eqref{case2-1} and \eqref{case2-2}.
 Now, the recursion \eqref{re} implies the following relations
 \begin{align*}
A_{m+1}(\beta_{m+1}+1-k)&=\sum_{i=0}^mA_m(\beta_m+2-k+i),\\[8pt]
A_{m+1}(\beta_{m+1}-k)&=\sum_{i=0}^mA_m(\beta_m+1-k+i).
 \end{align*}
Combining \eqref{indu4} and \eqref{ineq3}, we see that
 \begin{align}\label{g3}
\frac{\sum\limits_{i=0}^mA_m(\beta_m+2-k+i)}
{\sum\limits_{i=0}^mA_m(k-i)}
<\frac{\sum\limits_{i=0}^mA_m(\beta_m+1-k+i)}
{\sum\limits_{i=0}^mA_m(k+1-i)}.
 \end{align}
Taking into account altogether \eqref{g1},  \eqref{g2} and
\eqref{g3}, we conclude that
\begin{align}
\frac{A_{m+1}(\beta_{m+1}+1-k)}{A_{m+1}(k)}<
\frac{A_{m+1}(\beta_{m+1}-k)}{A_{m+1}(k+1)}, \ \ \ \ 2 \leq k \leq
r.
\end{align}
Plugging $k=r$ in above inequality gives
\[
\frac{A_{m+1}(r+2)}{A_{m+1}(r)}<\frac{A_{m+1}(r+1)}{A_{m+1}(r+1)}=1.
\]
Thus we have eventually established \eqref{want2}.  \qed

From the Theorem \ref{theorem}, we easily obtain the following
  log-concavity of $d_n(q)$:
\begin{cor}
The $q$-derangement number $d_n(q)$ $(n \geq 6)$ satisfies the
following property
\begin{align}\label{R_S}
\frac{A_n(1)}{A_n(2)}<\frac{A_n(2)}{A_n(3)}< \cdots
<\frac{A_n(r-2)}{A_n(r-1)}<\frac{A_n(r-1)}{A_n(r)},
\end{align}
where $r=\frac{n(n-1)}{2}-1$.
\end{cor}

 One should note that  the  sequence \eqref{R_S}
 does not involve the last term of $d_n(q)$ when $n$ is even.
However, if the last term is ignored, then $d_n(q)$ is log-concave
when  $n$ is even and $n\geq 6$. Moreover, it is clear that for $n$
even and $n\geq 6$ the last term of $d_n(q)$ must be excluded for
the consideration of log-concavity, because from Lemma \ref{lemma2}
 and the fact that $A_n(\beta_n)=1$, one sees that
 \[
\frac{A_n(\beta_n)}{A_n(\beta_n-1)}>\frac{A_n(\beta_n-1)}{A_n(\beta_n-2)},
 \]
 which violates the log-concavity.

\vspace{1cm}
 \noindent{\bf Acknowledgments.} We wish to thank Francesco Brenti
 and Michelle Wachs for their valuable comments.
 This work was supported by  the 973
Project, the PCSIRT Project of the Ministry of Education, the
Ministry of Science and Technology, and the National Science
Foundation of China.

\end{document}